\documentclass{llncs}
\usepackage{amsfonts}
\usepackage{amsmath}
\usepackage{mathtools}

\usepackage[cp1251]{inputenc}
\usepackage[english]{babel}

\usepackage{xspace}
\usepackage{graphicx}
\usepackage{cite}
\usepackage{hyperref}
\usepackage{color}

\newcommand{\conv}{\mathop{\rm conv}\nolimits}
\newcommand{\adj}{\mathop{\rm adj}\nolimits}

\newcommand{\RR}{{\mathbb R}}
\newcommand{\QQ}{{\mathbb Q}}

\newcommand{\ZZ}{{\mathbb Z}}
\newcommand{\BB}{{\mathcal B}}

\newcommand{\floor}[1]{\left\lfloor #1 \right\rfloor}

\newcommand{\fracc}[2]{\frac{\textstyle \strut #1 \strut}{\textstyle \strut #2 \strut}}

\newcommand{\FindIntpoints}{\text{{\sc FindIntpoints}\xspace}}
\newcommand{\ChooseNextNonIntegerVertex}{\text{{\sc ChooseNextNonIntegerVertex}\xspace}}
\newcommand{\GetCuts}{\text{{\sc GetCuts}\xspace}}

\newcommand{\diag}{\mathop{\rm diag}\nolimits}
\newcommand{\set}[1]{\left\{ #1\right\}}

\newenvironment{Code}{\medskip\frenchspacing\obeylines\obeyspaces}{\medskip}{\obeyspaces\global\let =\ \global\everymath{\let =\space}}

\newcommand{\While}{{\bf while{}}}

\newcommand{\If}{{\bf if{}}}

\newcommand{\Return}{{\bf return{}}}
\newcommand{\Procedure}{{\bf procedure{}}}
\newcommand{\For}{{\bf for{}}}

\begin{document}

\title{How to~Find the~Convex Hull of~All~Integer Points in~a~Polyhedron?\thanks{This work was supported by the Russian Science Foundation Grant No. 17-11-01336.}}

\titlerunning{How to~Find the~Convex Hull of~All~Integer Points?}  % abbreviated title (for running head)

\author{S.\,O.\,Semenov\orcidID{0000-0002-1498-3250} \and N.\,Yu.\,Zolotykh\orcidID{0000-0003-4542-9233}}

\authorrunning{S.\,O.\,Semenov \and N.\,Yu.\,Zolotykh} % abbreviated author list (for running head)

\institute{Mathematics of Future Technologies Center, \\ Lobachevsky State University of Nizhny Novgorod, Gagarin ave. 23, \\
Nizhny Novgorod 603950, Russia,\\
\email{sergey.semenov@itmm.unn.ru}, \email{nikolai.zolotykh@itmm.unn.ru}
}

\maketitle

\begin{abstract}
We propose a cut-based algorithm for finding all vertices and all facets of the convex hull of all integer points of a polyhedron defined by a system of linear inequalities. Our algorithm DDM Cuts is based on the Gomory cuts and the dynamic version of the double description method.
We describe the computer implementation of the algorithm and present the results of computational experiments comparing our algorithm with a naive one.

\end{abstract}

\section{Introduction}

In this paper we describe an algorithm DDM Cuts for finding all vertices and facets of the convex hull of all integer points of a convex polyhedron defined by a system of linear inequalities. This algorithm can be used as a subroutine in methods for solving convex and non-convex integer programming problems and for an experimental study of properties of combinatorial polyhedra.

Our method is based upon constructing Gomory cuts (see \cite{Hu,Schrijver,Shevchenko}). On each iteration the algorithm finds a non-integer vertex and constructs the cuts which cut off the vertex. The iterations continue until all vertices become integer. On each iteration we provide dual representation of the current polyhedron, i.e. we maintain both vertex and facet representation of the polyhedron. For this purpose we use a recently proposed dynamic variation of the double description method.

The double description method (DDM) is one of the most popular algorithms for finding the dual representation of a polyhedron \cite{MRTT1953,Chernikov1968en,Chernikova1965en}. The main difference between this algorithm and others for solving this problem is the use of both vertex and facet representations of the polyhedron on each iteration. DDM can quite naturally be extended to dynamic algorithms and we use one of such variations \cite{SemenovZolotykh2019} with additional heuristics proposed in \cite{Zolotykh2012en,BastrakovZolotykh2015en,BastrakovChurkinZolotykh2020}.

We describe a computer implementation of DDM Cuts and describe the results of computational experiment. The comparison of DDM Cuts with a naive algorithm (direct enumeration of all integer points and computing their convex hull) shows the superiority (in terms of performance) of DDM Cuts over the naive algorithm.

\section{Definitions}

Let $P$ be a {\em convex polyhedron} in $\RR^d$ given as $P=\set{x\in \RR^d:~ Ax\le b}$, where $A\in\RR^{m\times d}$ and $b\in\RR^m$. If $P$ is bounded then $P$ is a {\em convex polytope}, i.e. it can be represented as the convex hull of its vertices $v_1, v_2, \dots, v_n$.

Denote by $P_I$ the convex hull of all integer points contained in $P$, i.e. $P_I = \conv (P\cap \ZZ^d)$. We consider the problem of finding all vertices and all facets of $P_I$. Note that if $P$ is not bounded then $P_I$ can have infinitely many vertices and facets.

\section{Preliminaries}

In this section we consider the theoretical background of the proposed algorithm. This algorithm uses Gomory cuts.

An inequality $ax \le \beta$ is called the {\em cut} of a vertex $v$ of the polyhedron $P$ if $av > \beta$, but the inequality $ax \le \beta$ holds for all integer points of $P$. 

To construct all vertices of $P_I$ we can use the following procedure. First, find all vertices of $P$ using DDM or other algorithm.
If all vertices of $P_I$ are integer, then stop: they are all the vertices of $P_I$. Otherwise, take one of non-integer vertices and construct one or more cuts of this vertex. Join new inequalities to the original system and reconstruct the vertices of the polyhedron.

Let us consider the procedure for constructing Gomory cuts (see \cite{Hu,Schrijver,Shevchenko}).
Let $v$ be a non-integer vertex of $P=\set{x\in \RR^d:~ Ax\le b}$,
$A\in\ZZ^{m\times d}$, $b\in\ZZ^m$. Let $A_v\le b_b$ be a subsystem of 
$Ax\le b$, such that $A\in\ZZ^{d\times d}$, $\Delta = |\det A| \ne 0$ and $A_v v = b_v$. Consider the set
$$
\BB = \set{u\in\ZZ^d: ~ uA_v \equiv 0 \pmod{\Delta}} = \set{u=y \adj A_v:~ y\in \ZZ^d},
$$
where $\adj A_v = \det(A) A^{-1}$. It is known that $v\in\ZZ^d$
if and only if $ub_v \equiv 0 \pmod{\Delta}$ for any $u\in\BB$.
From this, it is not hard to obtain that if $u\in\BB$, $u\ge 0$ and $ub_v \not\equiv 0 \pmod{\Delta}$, then the inequality
$$
\frac{uA_v x}{\Delta} \le \floor{\frac{ub_v}{\Delta}}
$$
is a cut of $v$.
Moreover, in constructing the cuts it is enough to consider only vectors $u$ from the set $\overline{\BB}=\set{u\in\BB:~ 0\le u_j < \Delta~(j=1,2,\dots,d)}$.

Let $D$ be Smith normal form of $A_v$, i.e. $D=PA_vQ$, where
$P$ and $Q$ are unimodular matrices ($|\det P| = |\det Q| = 1$), 
$D=\diag(\delta_1, \delta_2, \dots, \delta_d)$ and $\delta_i$ is divided by $\delta_{i-1}$ $(i=2,3,\dots,d)$, then
$$
\overline{\BB} = \set{u = \left(\frac{\Delta}{\delta_1}\beta_1q_1 + \frac{\Delta}{\delta_2}\beta_2q_2 +  \dots + \frac{\Delta}{\delta_d}\beta_dq_d \right):~ 0\le \beta_j<\delta_j ~(j=1,2,\dots,d)},
$$
where $q_1, q_2,\dots,q_d$ are the columns of $Q$.
Moreover, by the vector $u\in\overline{\BB}$, the coefficients $\beta_1,\beta_2,\dots,\beta_d$ are determined uniquely.

\section{The proposed algorithm}

The proposed algorithm for finding the convex hull of all integer points in a polyhedron consists of iteratively cutting off non-integer vertices of the polyhedron while maintaining both its vertex and facet representations. 

First, we obtain the dual representation of the given polyhedron using DDM. Next, on each iteration, we choose a non-integer vertex and generate a set of cuts, that cut off the vertex from the polyhedron but keep all of its integer points. Those cuts are added to the facet representation of the polyhedron and its vertex description is updated accordingly via dynamic DDM \cite{SemenovZolotykh2019}. This process repeats until no non-integer vertices remain in the polyhedron.

Let $P$ be a polyhedron in $\RR^d$ given as $P = \{x \in \RR^d : Ax \leq b\}$, where $A \in \ZZ^{n \times d}$, $b \in \ZZ^n$. Then a procedure for finding the convex hull of integer points of $P$ can be described as follows.

%\begin{algorithmic}
%\Procedure{FindIntpoints}{$A$, $b$}
%\State Using DDM find $V \subset \QQ^d$, such that $\Conv V = \{x \in %\RR^d : Ax \leq b \}$
%\While{$\exists\; v \in V:$ $v \notin \ZZ^d $}
%	\State $v \in \QQ^d, A_v \in \ZZ^{d \times d}, b_v \in \ZZ^{d \times 1}$
%	\State $(A_v, b_v) \gets {\sc ChooseNextNonIntegerVertex}(V, A, b)$
%	\State $(A', b') \gets GetCuts(A_v, b_v) $
%	\State append rows of $A'$ to $A$
%	\State append rows of $b'$ to $b$
%	\State update $V$
%\EndWhile
%\State return $V$
%\EndProcedure
%\end{algorithmic}

\begin{Code}
\Procedure \FindIntpoints($A$, $b$):
        Using DDM find $V \subset \QQ^d$, such that $\conv V = \{x \in \RR^d :\, Ax \leq b \}$
        \While $\exists\; v \in V$, such that $v \notin \ZZ^d $:
                $v \in \QQ^d$, $A_v \in \ZZ^{d \times d}$, $b_v \in \ZZ^{d \times 1}$
                $A_v,\, b_v \gets \ChooseNextNonIntegerVertex(V,\, A,\, b)$
                $A',\, b' \gets \GetCuts(A_v,\, b_v) $
                append rows of $A'$ to $A$
                append rows of $b'$ to $b$
                update $V$
        \Return $V$
\end{Code}

The cuts added for each vertex $v$ are obtained via the procedure below, using the Smith normal form of the submatrix $A_v$ corresponding the vertex $v$, i.e. $A_v v = b_v$.

\begin{Code}
\Procedure \GetCuts($A_v$, $b_v$):
        $A' \gets \ZZ^{0 \times d}$
        $b' \gets \ZZ^0$
        Find Smith normal form $D$ of $A_v$, 
                i.e. $D = PA_{v}Q$, $P = (p_1, \dots , p_m)$, $D = \diag(d_1, \dots, d_m)$
        $\Delta \gets |\det(A_v)|$
        $C = \set{c = \sum\limits_{i=1}^m \fracc{\Delta}{d_i}\alpha_{i}p_i  \mod \Delta:~ \text{for all $\alpha_i$, such that } 0 \leq \alpha_i < d_i }$
        \For $c \in C$:
                \If $c \not\equiv 0 \pmod{\Delta}$ 
                        append $\fracc{cA_{v}x}{\Delta}$ to $A'$
                        append $\floor{\fracc{cb_v}{\Delta}}$ to $b'$
        \Return $A'$, $b'$
\end{Code}

For each vertex the algorithm generates and adds up to $|\det(A_v)|$ cuts. In order to minimize the number of cuts added on each iteration of the algorithm, we can use the following heuristics: we choose vertex $v$ with the minimum value of $\Delta$ (since vertices with large values of $\Delta$ might be cut off together with $v$).

\begin{Code}
\Procedure \ChooseNextNonIntegerVertex($V$, $A$, $b$):
        \For $v' \in V$:
            \If $v' \notin \ZZ^d$:
                    Find $A_{v'}$, $b_{v'}$, where $A_{v'}x\leq b_{v'}$ is a square subsystem
                            of $Ax\le b$, such that $A_{v'}v'= b_{v'}$
                    $\Delta' \gets |\det A_{v'}|$
                    \If $\Delta' < \Delta$:
                            $\Delta \gets \Delta'$, $A_v \gets A_{v'}$, $b_v \gets b_{v'}$
        \Return $(A_v, b_v)$
\end{Code}

The submatrix that represents vertex $v'$ is obtained by calculating the reduced row echelon form of the facets incident to $v'$. To further reduce the values of $\Delta$ on each iteration of the algorithm, each inequality (initial or added) is normalized by dividing all coefficients by their greatest common divisor.

\section{Computational Results} \label{sec_results}

A C++ implementation of the proposed algorithm has been developed. 
The computational experiments were performed on a computer with Intel(R) Core(TM) i7-8700K CPU at 3.70 GHz, Microsoft Windows 10 operating system, using the Microsoft Visual Studio 2017 compiler.

Also, a naive method for finding the convex hull of all integer solution to a system of linear inequalities was implemented. This algorithm directly enumerates all integer solutions and by means of DDM finds their convex hull.

For the experiment, random systems of linear inequalities were generated using the following procedure. Given $d$ and $k$, we generated
random uniformly distributed in the segment $[2k/3,\, k]$ integer values $a_1, a_2, \dots, a_d$. The algorithms were fed the system of linear inequalities
\begin{equation}
\left\{
\begin{array}{c}
a_1 x_1 + a_2 x_2 + \dots + a_d x_d \le k^2, \\[.45em]
x_j \ge 0 \qquad (j=1,2,\dots,d).
\end{array}
\right.
\label{f_exp_system}
\end{equation}
Given $d$, we took as $k$ each value in the logarithmically spaced sequence from $10$ to $1000$ of length $31$, i.e. $k= 10, 10\times 10^{1/15}, 10\times 10^{2/15}, \dots, 1000$. 

\begin{figure}[t]
    \centering
	\includegraphics*[width=0.85\linewidth]{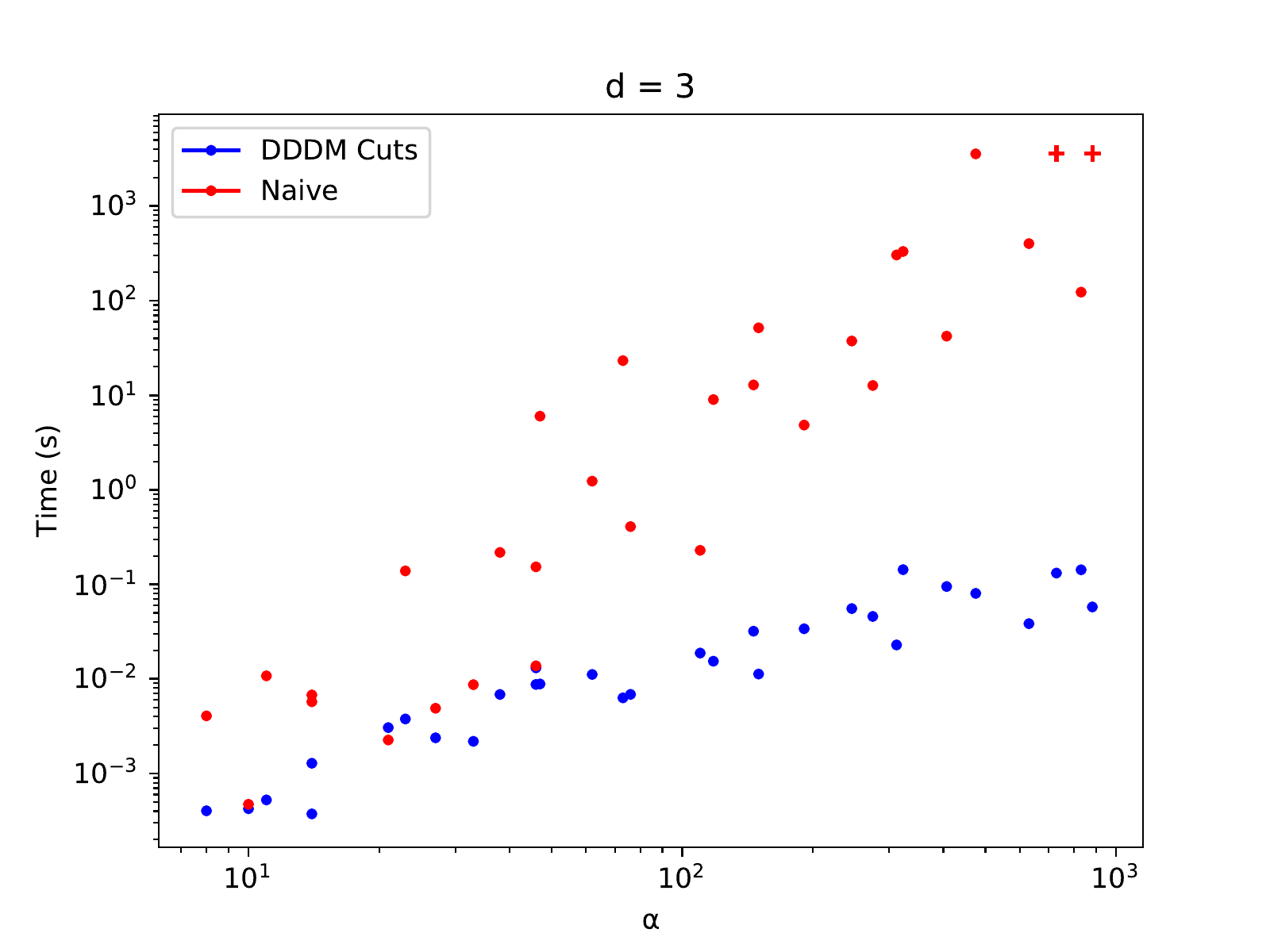}\\
	\caption{Performance of the proposed DDM Cuts algorithm and naive algorithm in solving the problem (\ref{f_exp_system}) for $d=3$. The dependence of the time (in seconds) on $\alpha=\max\set{a_1,a_2,\dots,a_d}$ is presented. Red $+$ sign means that the naive method was rejected due to exceeding the time limit of 600 sec.}
	\label{fig_time3}
\end{figure}

\begin{figure}
    \centering
	\includegraphics*[width=0.85\linewidth]{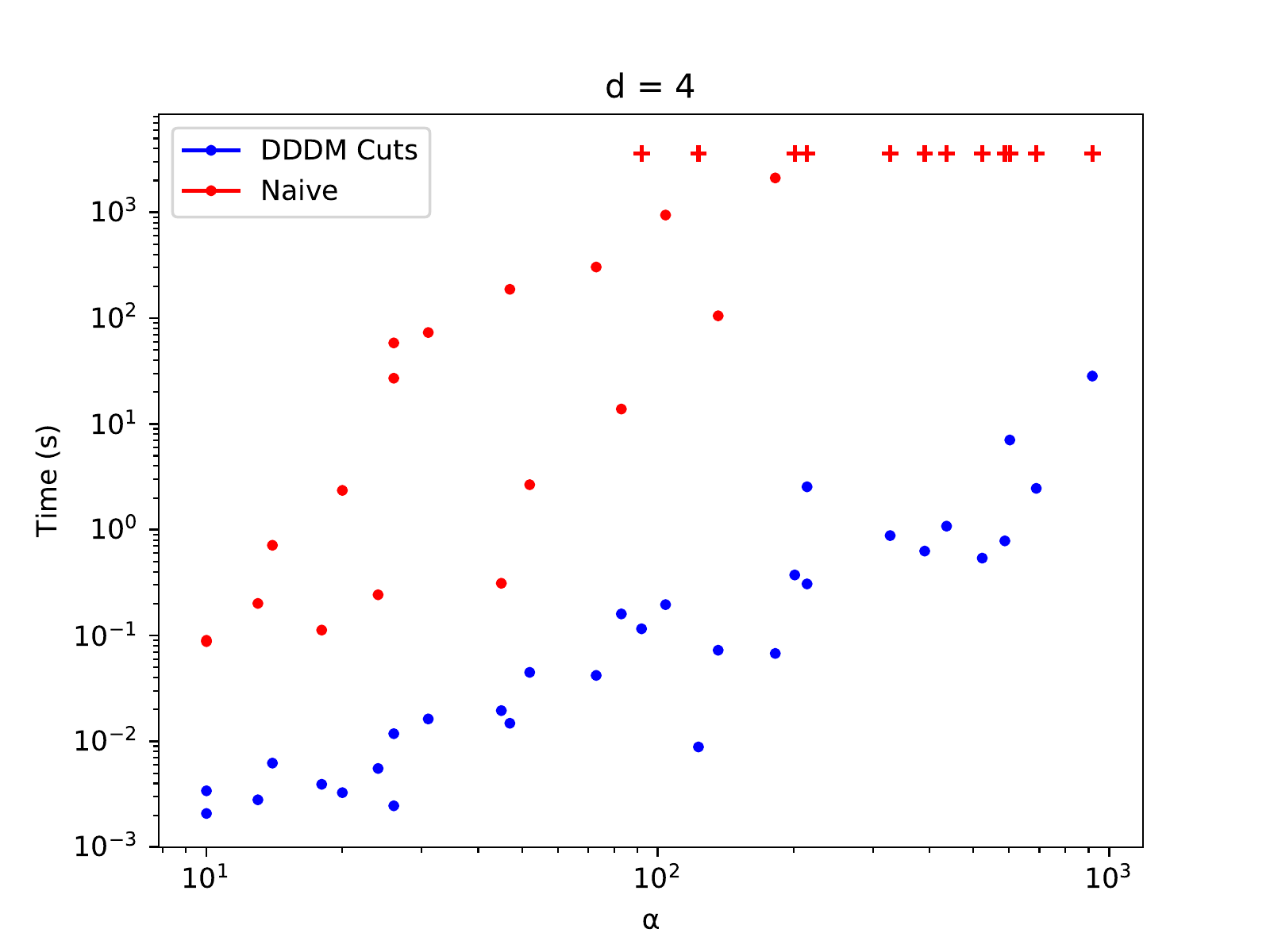}\\
	\caption{Performance of the proposed DDM Cuts algorithm and naive algorithm in solving the problem (\ref{f_exp_system}) for $d=4$. The dependence of the time on $\alpha$ is presented.}
	\label{fig_time4}
\end{figure}

The performance of the proposed DDM Cuts algorithm and naive algorithm for $d=3,4$ is presented on Figures~\ref{fig_time3},\,\ref{fig_time4}  respectively.
On the horizontal axis, we measure the value of $\alpha=\max\limits_j a_j$.
Apparently, the dependency of the execution time on $\alpha$ (if $d$ is fixed) is a power law.

\begin{figure}
    \centering
    \includegraphics*[width=0.85\linewidth]{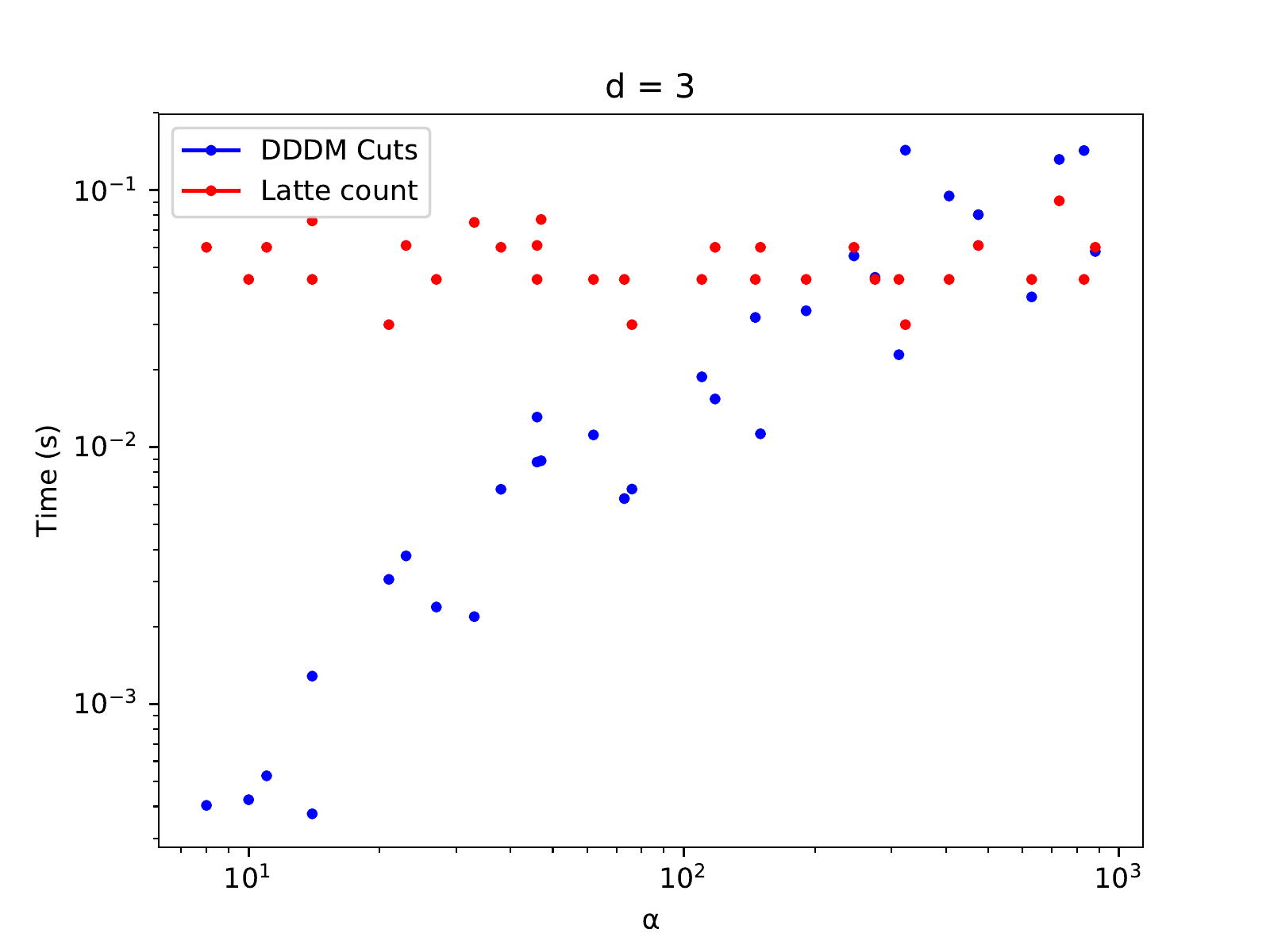}\\
    \caption{Performance of the proposed DDM Cuts algorithm and LattE in solving the problem (\ref{f_exp_system}) and counting of lattice points respectively for $d=3$. The dependence of the time on $\alpha$ is presented.}
    \label{fig_latte3}
\end{figure}

\begin{figure}
    \centering
    \includegraphics*[width=0.85\linewidth]{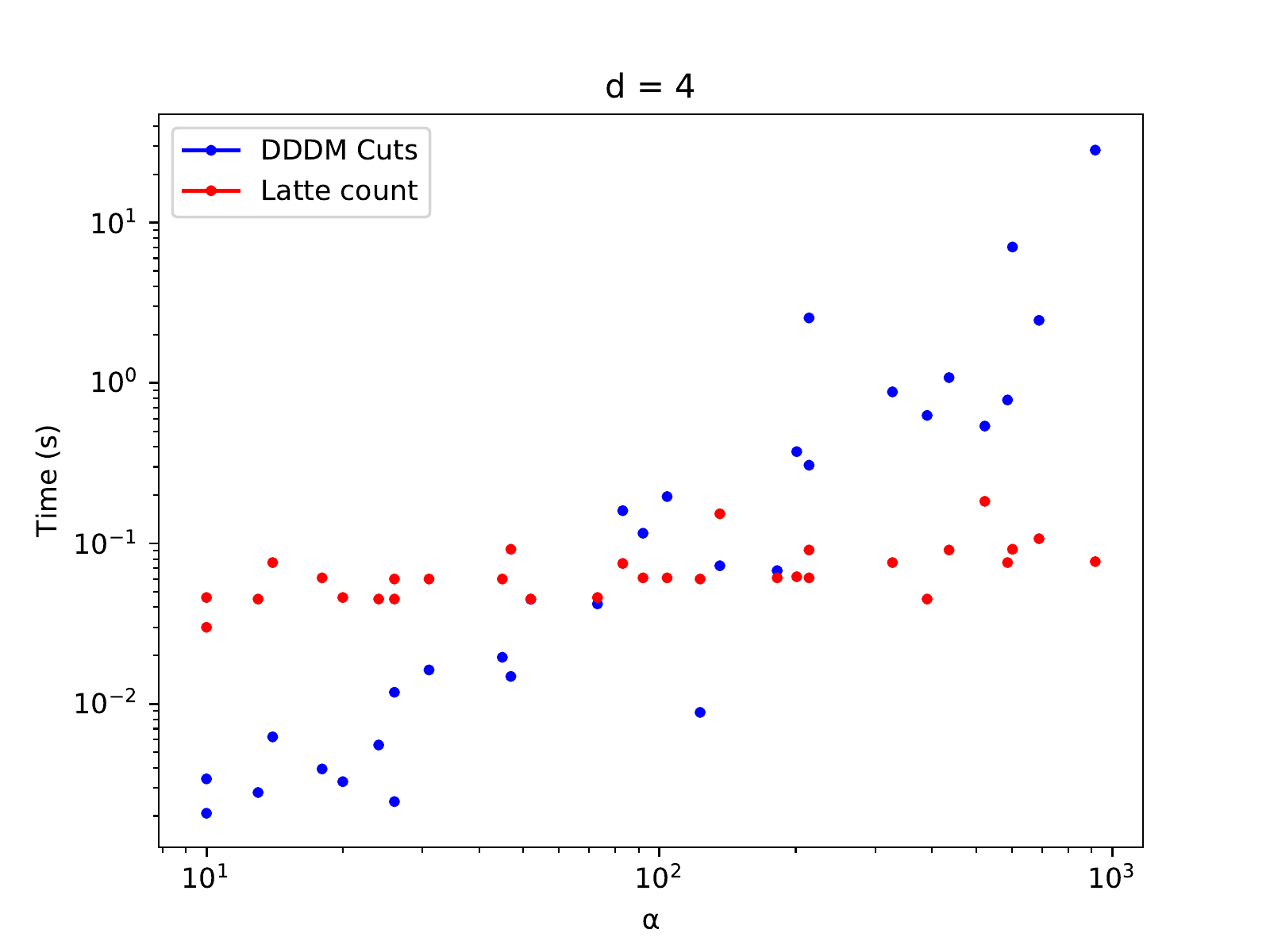}\\
    \caption{Performance of the proposed DDM Cuts algorithm and LattE in solving the problem (\ref{f_exp_system}) and counting of lattice points respectively for $d=4$. The dependence of the time on $\alpha$ is presented.}
    \label{fig_latte4}
\end{figure}

We're not aware of any other algorithms that tackle the same problem, so for the purposes of performance comparison we also consider an algorithm that solves a different one. LattE (Lattice points enumeration) \cite{DeLoeraHemmeckeTauzerYoshida} is an implementation of Barvinok's algorithm for counting lattice points inside convex polytopes. The performance of the proposed algorithm and LattE on the generated problem instances is presented on Figures~\ref{fig_latte3},\,\ref{fig_latte4}. LattE computation time is not dependent on the value of $\alpha$, so the proposed algorithm outperforms it only on lower values.

\begin{figure}
    \centering
	\includegraphics*[width=0.85\linewidth]{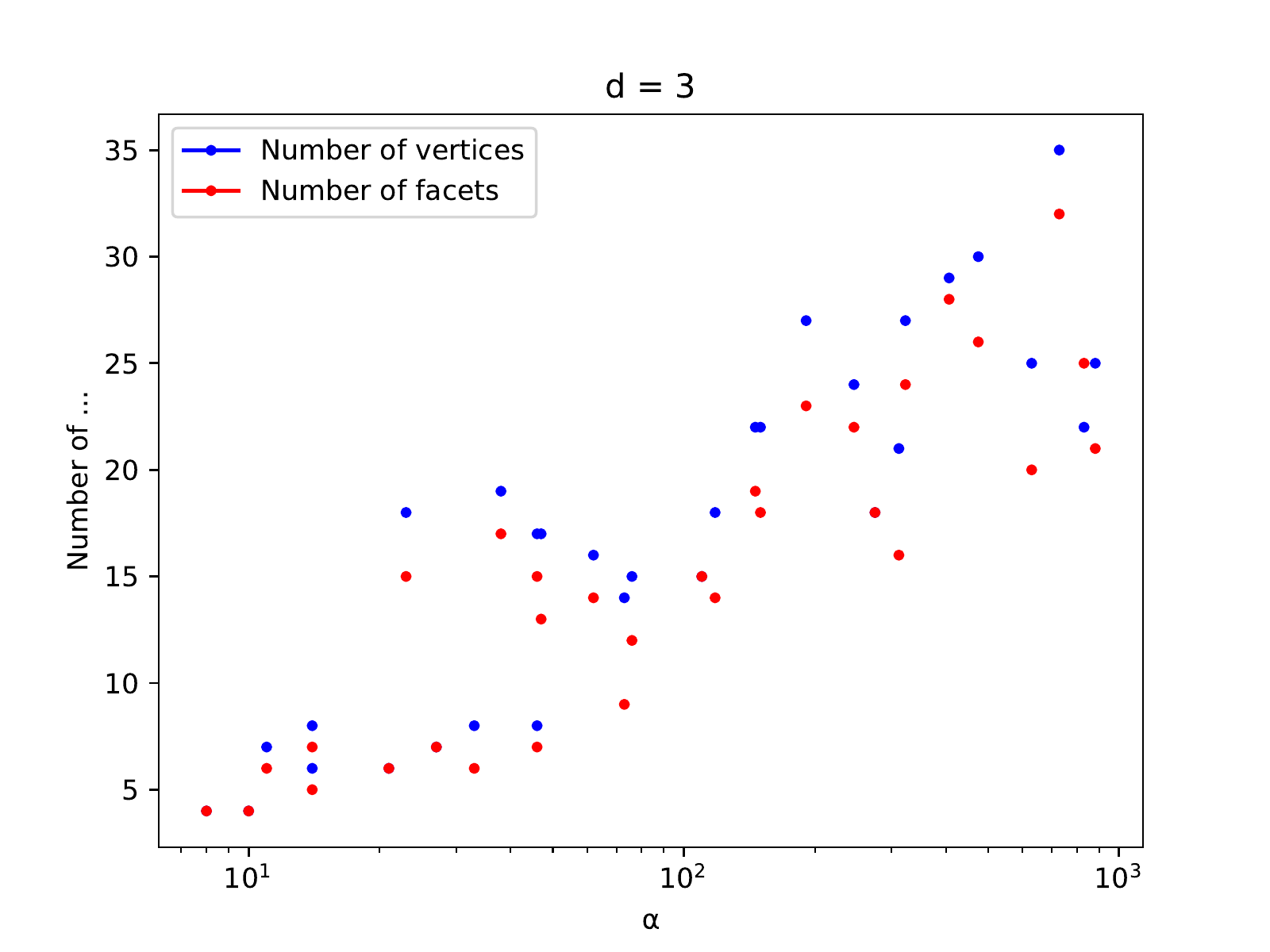}\\
	\caption{The number of vertices and facets for integer hull $P_I$ of all integer points of random polyhedron (\ref{f_exp_system}) for $d=3$.}
	\label{fig_vertexfacet3}
\end{figure}

\begin{figure}
    \centering
	\includegraphics*[width=0.85\linewidth]{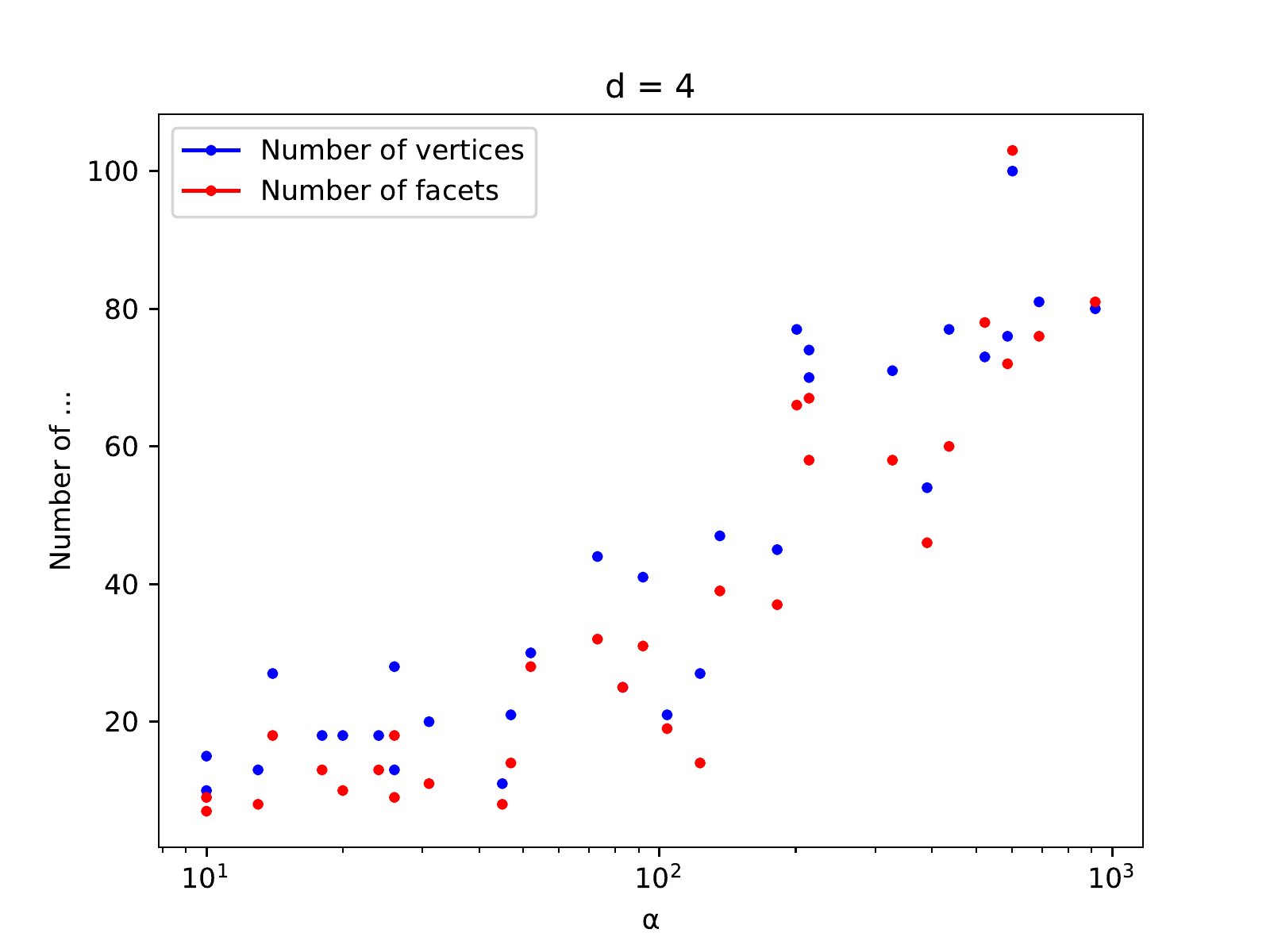}\\
	\caption{The number of vertices and facets for integer hull $P_I$ of all integer points of random polyhedron (\ref{f_exp_system}) for $d=4$.}
	\label{fig_vertexfacet4}
\end{figure}

The dependence of the number of vertices and facets of $P_I$ is represented on Figures~\ref{fig_vertexfacet3},\,\ref{fig_vertexfacet4}.

\begin{figure}
    \centering
	\includegraphics*[width=0.85\linewidth]{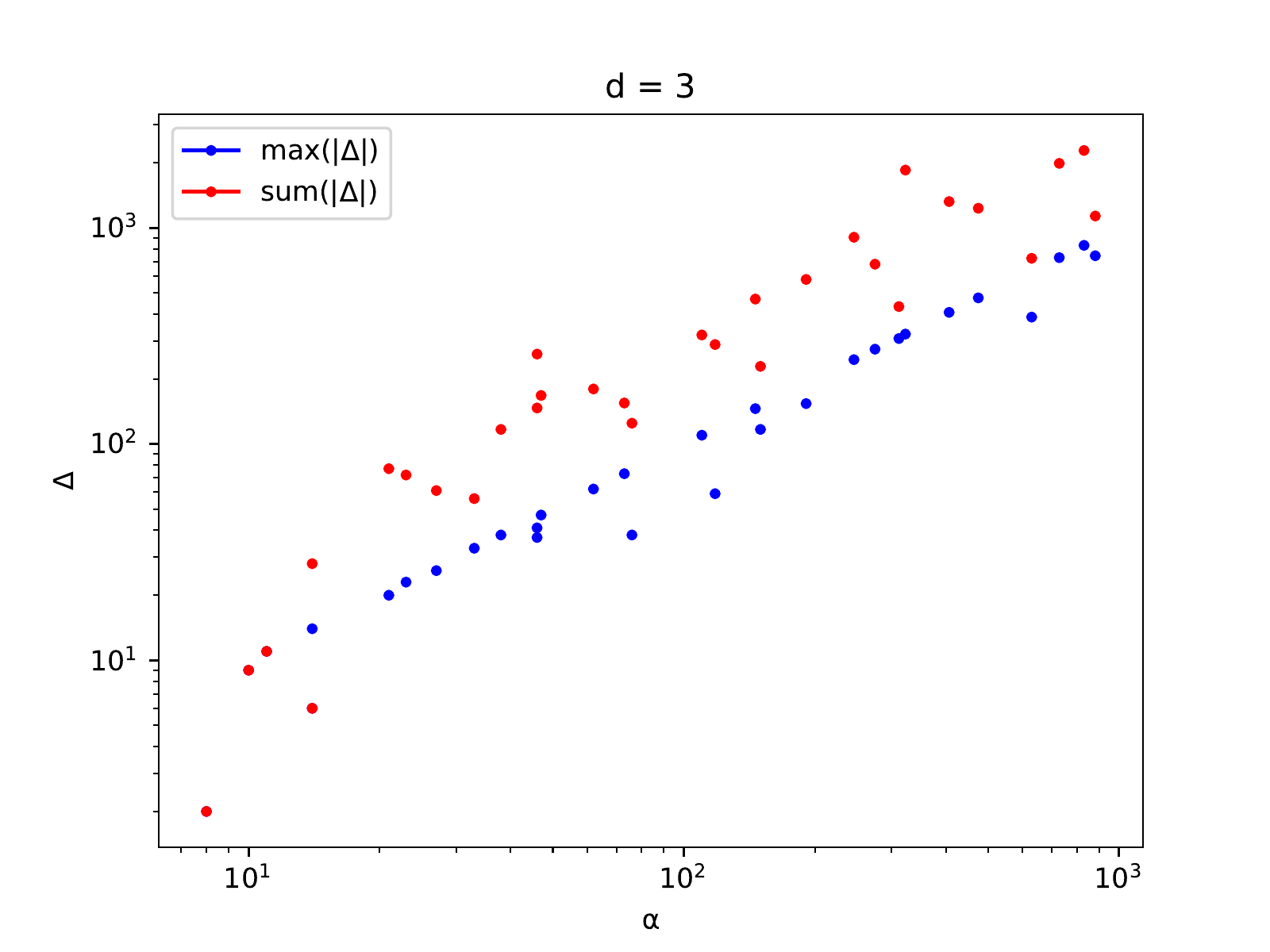}\\
	\caption{The maximal $\Delta$ and the sum of all $\Delta$ for all iterations of DDM Cuts algorithm, $d=3$.}
	\label{fig_det3}
\end{figure}

\begin{figure}
    \centering
	\includegraphics*[width=0.85\linewidth]{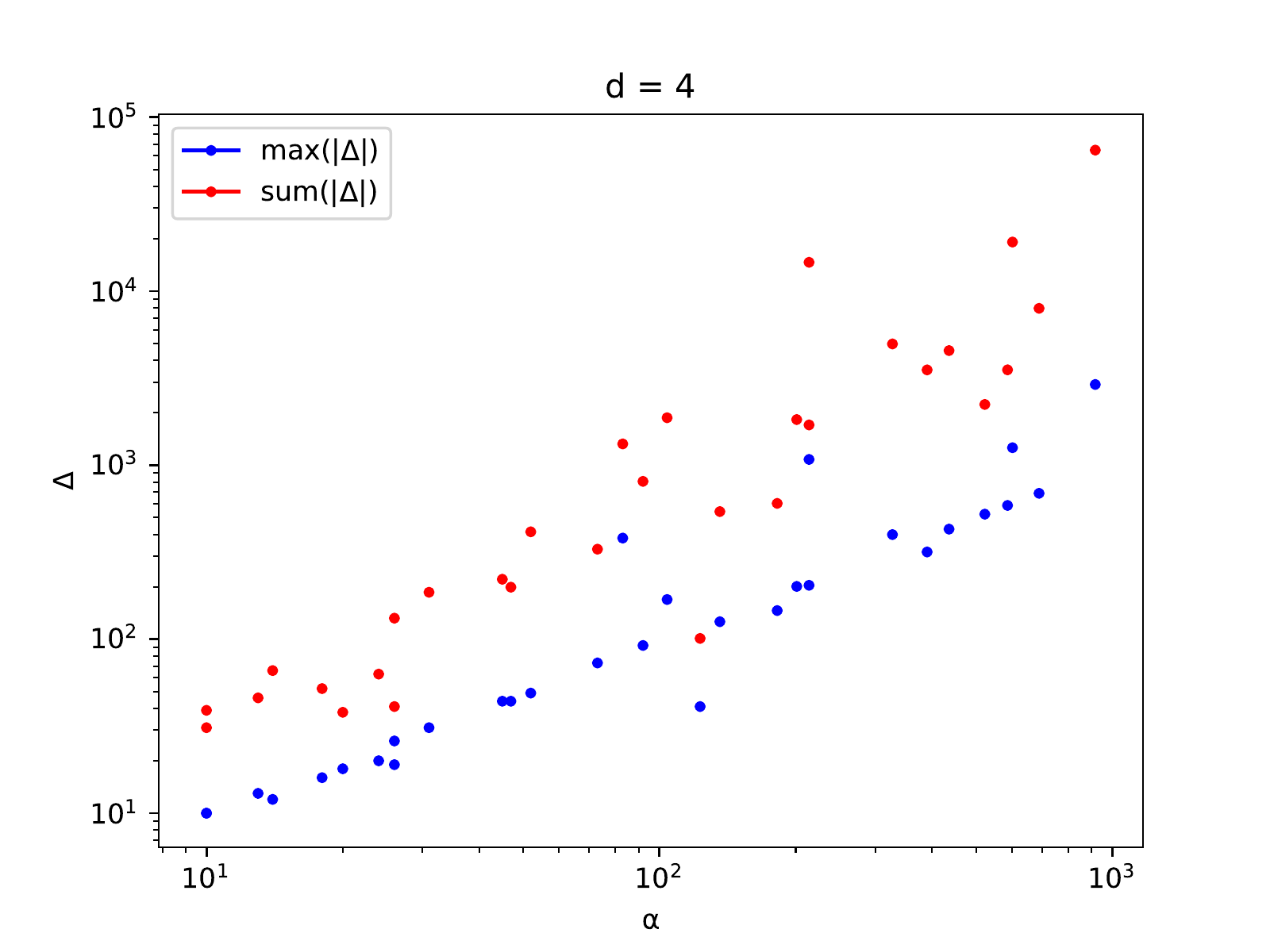}\\
	\caption{The maximal $\Delta$ and the sum of all $\Delta$ for all iterations of DDM Cuts algorithm, $d=4$.}
	\label{fig_det4}
\end{figure}

The dependence of the maximal $\Delta$ and the sum of all $\Delta$ for all iterations of DDM Cuts is represented on Figures~\ref{fig_det3},\,\ref{fig_det4}.

\section{Conclusion} \label{sec_conclusion}

In the paper we presented a cut-based algorithm for finding vertices and facets of the convex hull of all integer points of a polyhedron defined by a system of linear inequalities. The comparison in the framework of a computational experiment of our algorithm with a naive algorithm shows the superiority (in performance) of the former over the latter.

\end{document}